\documentclass{amsart}

\usepackage{graphicx}
\usepackage{amsmath}
\graphicspath{{./figures/}}

\usepackage{amssymb}
\usepackage{amsfonts}
\usepackage{amsmath}
\usepackage{amsthm}
\newtheorem{proposition}{Proposition} 

\newcommand\beq{\begin{equation}}
\newcommand\eeq{\end{equation}}

\renewcommand{\to}{\rightarrow}

\renewcommand{\div}{\operatorname{div}}

\newcommand{\tr}{\operatorname{tr}}

\newcommand{\ba}{\boldsymbol{a}} 
\newcommand{\bA}{\boldsymbol{A}} 

\newcommand{\be}{\boldsymbol{e}}

\newcommand{\f}{\boldsymbol{f}}
\newcommand{\bF}{\boldsymbol{F}}

\newcommand{\bG}{\boldsymbol{G}}

\newcommand{\bI}{\boldsymbol{I}} 

\def\bm{\boldsymbol m}

\newcommand{\bP}{\boldsymbol{P}}

\newcommand{\bu}{\boldsymbol{u}}

\newcommand{\bx}{\boldsymbol{x}} 

\newcommand{\bY}{\boldsymbol{Y}}

\newcommand{\bZ}{\boldsymbol{Z}}

\newcommand{\bphi}{\boldsymbol{\phi}}

\newcommand{\bsigma}{\boldsymbol{\sigma}}

\newcommand{\btau}{\boldsymbol{\tau}}

\def\Fcal{\mathcal{F}}
\def\Gcal{\mathcal{G}}

\newcommand{\NN}{\mathbb{N}}

 \newcommand{\R}{\mathbb{R}} 

\begin{document}

\title{Multi-material 
flow 
\\ with a new viscoelastic model} 

\author{Sébastien Boyaval}
\address{LHSV, EDF R\&D, ENPC, Institut Polytechnique de Paris, Chatou, France}
\address{Inria, France \hspace{1cm} email: sebastien.boyaval@enpc.fr}
%
\begin{abstract}
 We study a new quasilinear system of balance laws (i.e. conservative first-order PDEs with source terms) for the numerical simulation by a front-capturing approach of \emph{viscoelastic} athermal multi-material flows with sharp interfaces (like a mixture of immiscible fluids and solids). It naturally extends with an energy balance that is meaningful to stability and thermodynamics. In particular, the unequivocal definition of solutions to Cauchy problems can be guaranteed in some cases, thanks to a symmetric-hyperbolic formulation.
\end{abstract}

\maketitle

\section{Introduction}

The numerical simulation of fast 
mechanical waves throughout \emph{a heterogenenous continuum with variable wave speeds} is an active research topic.
It is in particular true in the case with sharp interfaces between the heterogeneities,
like the propagation of a shock thoughout a mixture of finitely-many 
immiscible materials\footnote{We call here ``material'' what is sometimes called a ``phase'', and which coincides here with a subset of the simulation domain with one specific mechanics. 
} 
see e.g. \cite{MR3481008,Rodriguez2019,MR4038165}.

One 
modelling approach for multi-material flows is with \emph{interface-capturing methods}, when the mechanics of all materials are evolved simultaneously with a single system of bulk equations for the whole simulation domain. 
In the interface-capturing approach, one typically expects that the interfaces between the immiscible materials are well-defined implicitly, as discontinuities in the materials properties, and that they at worst require only purely-numerical reconstruction at each time-step after discretization (i.e. no additional physics). As long as the essential multi-fluid physics can be ensured through one single system of bulk equations, 
the latter methods are an appealing alternative to the 
methods that simulate multi-materials flows with one system of equations per material subdomain and that \emph{explicitly} propagate sharp interfaces using coupling conditions between 
different materials (therefore as numerous as the variety of interfaces). 
However, various systems of bulk equations 
have been used by interface-capturing methods to simulate multi-material motions on times $t \in (0,T)$ given a force field $\f(t,\bx)$ in the ambiant physical space $\bx\in\R^3$, see e.g. \cite{MR3481008,Rodriguez2019,MR4038165}, and it is not clear yet what is best to handle correctly a given configuration in general. 

Here, we propose a new 
system of PDEs to govern athermal 
\emph{viscoelastic} multi-material flows of immiscible fluid-solid mixtures with sharp interfaces, therefore \emph{possibly with dissipation} 
; it is a quasilinear system of balance laws (i.e. conservative first-order PDEs with source terms) such that the stability of solutions to Cauchy problems can be studied starting with the question of
hyperbolicity \cite{DafermosBook4}.
%

In the sequel, we first introduce our new system of quasilinear PDEs (\ref{eq:momentum_spatial_tensor}--\ref{gammalaw}--\ref{eq:reference_spatial_tensor}--\ref{eq:K_spatial_tensor}--\ref{formula}--\ref{eq:deformation_spatial_tensor}--\ref{eq:G_spatial_tensor}--\ref{eq:lambda_spatial_tensor}--\ref{eq:Z_spatial_tensor})
to model multi-material flows (fluid-solid mixtures) with sharp interfaces and possibly viscosity.
Then we discuss the unequivocal definition of solutions to Cauchy problems through the standard approach relying on a hyperbolic formulation.
Finally, some conclusions can be drawn about the stability of solutions justifying interface-capturing simulations. In the case without viscosity, our interface-capturing approach seems to be the first rigorous one that guarantees the unequivocal definition of solutions to Cauchy problems.

\section{The new viscoelastic system} 

In the literature which we are aware of, the interface-capturing methods 
complement the 
general momentum balance governing the velocity $\bu(t,\bx)$ 
\begin{equation}
\label{eq:momentum_spatial_tensor}
\rho( \partial_t \bu +  \bu\cdot\nabla \bu ) - \nabla\cdot\bsigma = \rho \f
\end{equation}
with various definitions of the mass density $\rho(t,\bx)$ and of the Cauchy stress $\bsigma(t,\bx)$ for almost all $t$ and $\bx$ (in the Lebesgue measure sense).
Hypoelastic closures use $\bsigma$ as a variable, see for instance 
\cite{Rodriguez2019} (which encompasses our scope, insofar as it considers miscible cases with 
fraction variables $\alpha\in(0,1)$ entering the definition of $\bsigma$). 
Here, we follow \cite{MR3481008,MR4038165} and consider only 
hyperelastic formulations. 

Then, in a ``perfect fluid'' subdomain for instance, $\bsigma$ should coincide with a spherical stress $-p(J)\bI$
where the pressure $p=-\partial_J e_K$ is derived from an 
internal energy function of the volumetric deformation $J:=\frac{\rho_r}\rho$ with respect to a reference mass density $\rho_r$ of the material like 
\begin{equation}
\label{gammalaw} p=K J^{-\gamma}
\end{equation}
when $K$ is the reference value of the pressure. 
Other pressure formulas are possible: Van der Waals law e.g. 
see \cite{MR3481008} or $p=K (J^{-\gamma}-J)+p_r$ so $p=p_r$ when $J=1$ whatever $K$. 
In any case, \eqref{eq:momentum_spatial_tensor} is 
closed for a perfect fluid 
after complementing the definition $\bsigma=-p\bI$ and 
\eqref{gammalaw} by
\begin{equation}
\label{eq:dilatation_spatial_tensor} \partial_t J^{-1} + \nabla\cdot( \bu J^{-1} ) = 0
\end{equation}
to govern the volume of matter, 
plus 
\begin{align}
\label{eq:reference_spatial_tensor} \partial_t \rho_r + (\bu\cdot\nabla) \rho_r & = 0
\\
\label{eq:K_spatial_tensor} \partial_t K + (\bu\cdot\nabla) K & = 0
\end{align} 
to define the material properties $\rho_r$ 
and $K$, 
possibly with fictitious values outside the perfect-fluid subdomain (like in ``ghost-fluid'' methods) if the latter variables are not relevant elsewhere.

Here, we propose that Cauchy stress is everywhere of \emph{neo-Hookean/Mooney-Rivlin type}\footnote{Note that our stress formula does not assume a strict 
separation of the volumetric and isochoric contributions, as opposed to most neo-Hookean/Mooney-Rivlin formulations.} following \cite{MR3481008}, and we specifically propose the following \emph{new} formula
\begin{equation}
\label{formula} \bsigma = 
-p\bI + \btau\,, \btau:=J^{-1}(\bF\bA\bF^T-G\bI) 
\end{equation}
where $\bF$ is the usual deformation gradient with respect to a reference state 
thus governed by
\begin{equation}
 \label{eq:deformation_spatial_tensor} (\partial_t + \bu\cdot\nabla) \bF = (\nabla\bu)\bF
\end{equation}
and such that its determinant $|F|\equiv J$ satisfies \eqref{eq:dilatation_spatial_tensor}, while the material properties are defined using additional variables and additional equations like
\begin{align}
\label{eq:G_spatial_tensor} \partial_t G + (\bu\cdot\nabla) G & = 0
\\
\label{eq:lambda_spatial_tensor} \partial_t \lambda + (\bu\cdot\nabla) \lambda & = 0
\\
\label{eq:Z_spatial_tensor} \partial_t \bA + (\bu\cdot\nabla) \bA & = \frac{G\bF^{-1}\bF^{-T}-\bA}\lambda
\end{align}
as long as the matrix-valued variable $\bA$ is symmetric-postive.
Note that the tensor variable $\bA$ 
in \eqref{formula} is a contribution to the \emph{second Piola-Kirchhoff stress} 
$J\bF^{-1}\bsigma\bF^{-T}$ that is related to the strain tensor $\bF^{-1}\bF^{-T}$ 
and a variable elastic-modulus $G$ through the relaxation 
\eqref{eq:Z_spatial_tensor} with a charateristic material time $\lambda$, so that the resulting extra-stress $\btau$ in \eqref{formula} can be interpreted as viscoelastic insofar as it satisfies a constitutive law of Maxwell-type
\begin{equation}
\label{eq:extrastress_spatial_tensor}
(\partial_t + \bu\cdot\nabla) \btau - \nabla\bu\btau - \btau\nabla\bu^T + \btau\div\bu = \frac{\mu(\nabla\bu + \nabla\bu^T)-\btau}\lambda
\end{equation}
with a viscosity $\mu=J^{-1}G\lambda$ -- this is why we call our 
approach \emph{viscoelastic}. We shall discuss later the choice of the exponent $-\tfrac1{r_Z}$.

The system of quasilinear PDEs (\ref{eq:momentum_spatial_tensor}--\ref{gammalaw}--\ref{eq:reference_spatial_tensor}--\ref{eq:K_spatial_tensor}--\ref{formula}--\ref{eq:deformation_spatial_tensor}--\ref{eq:G_spatial_tensor}--\ref{eq:lambda_spatial_tensor}--\ref{eq:Z_spatial_tensor}) 
is directly inspired by our previous works \cite{Boyaval2021,boyaval2024,S0219891622500096,LEE2026118742} about a symmetric-hyperbolic system of viscoelastic balance laws for (single-material) flows with homogeneous material properties. It is a Eulerian description of a viscoelastic continuum generalized to \emph{multi-material} flows, when the material properties $\rho_r,K,G,\lambda$ 
vary in addition to $\bA$ (already defined as a variable in \cite{Boyaval2021,boyaval2024,S0219891622500096,LEE2026118742} to introduce viscous dissipation through a relaxation equation).
It is relevant to \emph{immiscible} fluid-solid mixtures insofar as $\rho_r,K,G,\lambda$ 
are purely advected following pure transport equations (\ref{eq:reference_spatial_tensor}--\ref{eq:K_spatial_tensor}--\ref{eq:G_spatial_tensor}--\ref{eq:lambda_spatial_tensor}).
Moreover, as long as its solutions are well-defined, 
our hyperelastic 
interface-capturing approach \emph{actually defines multi-material flows at the PDE level whatever the (number of) values assumed by the material properties} $\rho_r,K,G,\lambda$ 
functions of $t$ and $\bx$ (within a reasonable regularity class like BV functions), without the help of additionnal hypotheses 
by contrast with other approaches for immiscible fluid-solid mixtures \cite{MR3481008,MR4038165}.
In particular, it \emph{does not need the introduction of additional geometrical variables to track material interfaces} (like interface-tracking methods), as long as interfaces are well-defined 
by unequivocal (jumps of) piecewise-smooth solutions are well-defined (e.g. in 1D hyperbolic cases),
as one can expect from interface-capturing methods as opposed to interface-tracking methods.
So, as regards the numerical computation of a solution, with a view to simulating immiscible fluid-solid mixtures neglecting ``interfacial stresses'' 
but not necessarily viscous stresses, the effort required by our interface-capturing approach is independent of (i.e does not scale with) the number of interfaces ``types''.

Let us now study the unequivocal definition of 
solutions to the Cauchy problem 
with (\ref{eq:momentum_spatial_tensor}--\dots--\ref{eq:Z_spatial_tensor}) such that both $\bF$ \emph{and} $J$ are (independent) variables, and the the Piola idendities
\begin{equation}
\label{piola} \div(J^{-1}\bF^T) = 0                                                                                                             \end{equation}
are satisfied. Indeed, to study the unequivocal definitions of solutions to the quasilinear PDEs (\ref{eq:momentum_spatial_tensor}--\dots--\ref{eq:Z_spatial_tensor}), here we shall concentrate on establising a symmetric-hyperbolic formulation and we therefore consider from the start our visco(hyper)elastic multi-material flows in the standard ``polyconvex'' framework where the hyperbolicity of standard (single-material) elastodynamics can be studied \cite{DafermosBook4,wagner-2009}.

\section{Hyperbolicity}
\label{sec:hyperbolic}

First, note that smooth solutions to (\ref{eq:momentum_spatial_tensor}--\dots--\ref{eq:Z_spatial_tensor}) are equivalently solutions to a quasilinear system of \emph{balance laws}, i.e. conservative first-order PDEs plus algebraic source terms,  as long as \eqref{piola} is satisfied and $J^{-1}\bF$ coincides with the 
transposed adjugate of $\bF^{-1}$ i.e. one considers
\begin{equation}
\label{eq:deformation_material_tensor_cons}
\partial_t (J^{-1}\bF) + \div\left(\bu\otimes J^{-1}\bF-J^{-1}\bF\otimes\bu\right)=0
\end{equation}
instead of \eqref{eq:deformation_material_tensor} see e.g. \cite{DafermosBook4,wagner-2009,Boyaval2021}.
The Piola idendities \eqref{piola} are preserved by smooth solutions to \eqref{eq:deformation_material_tensor_cons}, and using \eqref{eq:dilatation_spatial_tensor}, balance laws can be formulated like 
$$
\partial_t \rho + \div\left(\rho\bu\right)=0
$$
for $\rho:=\rho_r J^{-1}$ on recalling \eqref{eq:reference_spatial_tensor}, or
$$
\partial_t \left( \rho_r^{r_\rho} J^{-1} \right) + \div\left(\rho_r^{r_\rho} J^{-1}\bu\right)=0 \,.
$$
We shall discuss later adequate choices for $r_\rho$ as well as for $r_K,r_G,r_\lambda$  
in
\begin{align}
\partial_t \left( K^{r_K} J^{-1} \right) + \div\left( K^{r_K} J^{-1}\bu\right) &=0
\\
\partial_t \left( G^{r_G} J^{-1} \right) + \div\left( G^{r_G} J^{-1}\bu\right) &=0
\\
\partial_t \left( \lambda^{r_\lambda} J^{-1} \right) + \div\left( \lambda^{r_\lambda} J^{-1}\bu\right) &=0
\end{align}
but in any case, the following conservation of linear momentum also holds
\begin{equation}
\label{eq:momentum_law}
\partial_t ( \rho \bu ) +  \nabla ( \rho \bu\otimes\bu ) - \nabla\cdot\bsigma = \rho \f
\end{equation}
and we therefore introduce $\bm:=\rho_r\bu$ so $\rho\bu=J^{-1}\bm$.

Now, a system of balance laws like that satisfied by
\begin{equation}\label{variable}
(J^{-1},J^{-1}\bm,J^{-1}\bF,J^{-1}\rho_r^{r_\rho},J^{-1}K^{r_K},J^{-1}G^{r_G},J^{-1}\lambda^{r_\lambda},J^{-1}\bA^{r_A}) 
\end{equation}
has a symmetric-hyperbolic formulation if 
it 
admits a strictly convex extension i.e. if
\begin{equation}
\label{extension} \partial_t E + \div\bG = H                                                                                                             \end{equation}
is a balance law satisfied by a strictly convex functional $E$ of the 
variable \eqref{variable} with smooth algebraic functionals $\bG$, $H$ of 
\eqref{variable}  \cite{DafermosBook4,wagner-2009}.
Furthermore, our hyperelastic modelling approach 
naturally yields a non-trivial additional balance law \eqref{extension} for the system satisfied by \eqref{variable}, with
\begin{equation} \label{energy}
E := J^{-1} \left( \frac{|\bm|^2}{2\rho_r} + K \frac{J^{1-\gamma}}{\gamma-1} 
+ \frac{\bF^T\bF:\bA}2 - G \log J  - \frac{G}2 \log|A|  \right)
\end{equation}
$\bG = E\bu -\bsigma\bu$ and an algebraic source term $H\le J^{-1}\bm\cdot\f$.
Note that \eqref{extension} is a formulation of the second principle of thermodynamics for our new system modelling athermal mechanics\footnote{
 In particular, it is the requirement of such a second principle of thermodynamics that limits the validity of our formulation to \emph{positive-definite} symmetric $\bA$, so as to define \eqref{extension} properly, and then $\bA^{r_A}$ as well.
}, thanks to the sign $H\le0$ when $\f=0$.
So the system of balance laws for \eqref{variable} is symmetric-hyperbolic when $E$ is strictly convex with respect to the conservative variable \eqref{variable}.
However, to identity a region of (strict) convexity for $E$ with respect to a set of conserved variables,
a number of adjustments are still in order, notably the choice of conserved variables. 

The determinant $|A|$ in 
\eqref{energy} is neither convex nor concave with respect to $\bA$. So, to handle the convexity of $E$, we 
shall add an independent variable like
$J^{-1}|A|$ to \eqref{variable} and 
\begin{equation}
\label{eq:Zdet_spatial_tensor} 
(\partial_t + \bu\cdot\nabla) |A| =
|A|\frac{G\bF^{-1}\bF^{-T}:\bA^{-1}-d}\lambda
\end{equation}
to the system of balance laws, in the ``polyconvexity'' spirit, $d=\tr\bI=3$ being the space dimension.
Let us also recall that for smooth flows such that it holds $J\equiv|F|$, \eqref{piola} and
$$ \bu(t,\bx) = \partial_t \bphi_t(\bphi_t^{-1}(\bx)) \quad \bF(t,\bx) = \nabla_{\ba} \bphi_t(\bphi_t^{-1}(\bx)) $$
with a diffeomorphism $\bphi_t^{-1}:\bx\in\R^3 \to \ba\in\R^3$ from the ambiant physical space to a reference ``material'' configuration space, simpler balance laws can be derived
\begin{align}
\label{eq:momentum_material_tensor} \partial_t \bm - \nabla_{\ba}\cdot\bP & = \rho_r \f
\\
\label{eq:dilatation_material_tensor} \partial_t J + \nabla_{\ba}\cdot\left( \hat\bF^T \bm \rho_r^{-1} 
\right) & = 0
\\
\label{eq:deformation_material_tensor} \partial_t \bF - \nabla_{\ba} \left( \bm \rho_r^{-1} 
\right) & = 0
\\
\label{eq:reference_material_tensor} \partial_t \rho_r^{r_\rho} & = 0
\\
\label{eq:K_material_tensor} \partial_t K^{r_K} & = 0
\\
\label{eq:G_material_tensor} \partial_t G^{r_G} & = 0
\\
\label{eq:lambda_material_tensor} \partial_t \lambda^{r_\lambda} & = 0
\end{align}
using material coordinates. When $r_A$ is any power, source terms is complex in general e.g.
\begin{equation}
\label{eq:Z_material_tensor_rA} \partial_t \bA^{r_A} = \left(\frac1{2i\pi}\int_{\Gamma}s^r(s\bI-\bA)^{-1}\bF^{-1}\bF^{-T}(s\bI-\bA)^{-1}ds-\bA^{r_A}\right) \lambda^{-1}
\end{equation}
using a contour $\Gamma$ enclosing the eigenvalues of $\bA$ in the complex plane, but choosing $r_A$ as a negative integer will turn out natural and will therefore simplify \eqref{eq:Z_material_tensor_rA} in the sequel.
As usual, we have abusively denoted $\bu$ above the new variable $\bu(t,\bphi_t(\ba))$ mapped from $\bu(t,\bx)$ insofar as the subscript $_{\ba}$ 
avoids misunderstanding, and in \eqref{eq:momentum_material_tensor} 
\begin{equation}
\label{piolakirchhoff} \bP = 
\left( - K J^{-\gamma} - G J^{-1}\right)\hat\bF + \bF\bA
\end{equation}
denotes the first Piola-Kirchhoff tensor as usual.
The latter system (\ref{eq:momentum_material_tensor}--\dots--\ref{eq:Z_material_tensor_rA}) defines the so-called \emph{Lagrangian description} equivalent to the Eulerian description (\ref{eq:momentum_spatial_tensor}--\dots--\ref{eq:Z_spatial_tensor}) insofar as their solutions are in one-to-one correspondance through $\phi_t$.
Note that, the constraint to be imposed on initial conditions and preserved by the flow, i.e. \eqref{piola} in the Eulerian description, becomes
\begin{equation}
\label{piola_material} \nabla_{\ba}\cdot\hat\bF^T = 0
\end{equation}
in the Lagrangian description ; it is essential to the conservative formulation \eqref{eq:dilatation_material_tensor} for $J\equiv|F|$.
The extension \eqref{extension} 
with $\bG = E\bu -\bsigma\bu$ and an algebraic source term $H\le J^{-1}\bm\cdot\f$ becomes
\begin{equation}
\label{extensionbis} \partial_t (JE) - \div_{\ba}(\bu\cdot\bP) = JH                                                                                                             \end{equation} 
in material coordinates, and indeed the convexity of $E$ with respect to Eulerian (conservative) variables 
is equivalent to the convexity of $JE$ with respect to their Lagrangian counterparts 
i.e. hyperbolicity is preserved by the smooth change of variable and coordinates above see e.g. \cite{wagner-2009}.

Now, recall $(x,y)\to x^\alpha y^\beta$ is convex strictly
\begin{itemize}
 \item on $\R_{>0}^2\ni(x,y)$ if $\alpha,\beta<0$, and
 \item on $\R\times\R_{>0}\ni(x,y)$ if $\beta<0$, $\alpha>1-\beta>1$. 
\end{itemize}
So $(\rho_r^{r_\rho},\bm)\to\frac{|\bm|^2}{2\rho_r}$ is strictly convex on $\R_{>0}\times\R^d$ 
whatever $r_\rho>0$, and $(K^{r_K},J)\to K \frac{J^{1-\gamma}}{\gamma-1}$ is strictly convex on $\R_{>0}\times\R_{>0}$ when $\gamma>1$ and $r_K\in(0,\tfrac1{\gamma})$.
Moreover, we recall $(\bF,\bZ)\to\bF^T\bF:\bZ^{-\tfrac1{r_Z}}$ is strictly convex on $\R^{d\times d}\times \R^{d\times d}_{>0}$ when $r_Z>1$, $\R^{d\times d}_{>0}$ denoting the set of positive symmetric real-valued matrices \cite{lieb-1973,Boyaval2021}. 
We therefore claim at this stage:
\begin{proposition}\label{prop:perfect}
Assume $\gamma>1$ and $r_K\in(0,\tfrac1{\gamma})$, $r_\rho>0$, $r_A<-1$.\\
When $\lambda\to\infty$, $G\equiv\frac{J \mu}{\lambda}\to$, the solutions to the system of quasilinear balance laws (\ref{eq:momentum_material_tensor}--\ref{eq:dilatation_material_tensor}--\ref{eq:deformation_material_tensor}--\ref{eq:reference_material_tensor}--\ref{eq:K_material_tensor}--\ref{eq:Z_material_tensor_rA}) that satisfy \eqref{piola_material} also satisfy a strictly convex extension on
\begin{equation}\label{variable2}
\R_{>0}\times\R^d\times\R^{d\times d}\times\R_{>0},\R_{>0},\R_{>0}^{d\times d} \ni
(J,\bm,\bF,\rho_r^{r_\rho},K^{r_K},\bA^{r_A}) 
\end{equation}
and therefore define an unequivocal Lagrangian description for \emph{perfect} flows of immiscible infinitely-many materials insofar as they are the solutions to well-defined Cauchy problems on $(0,T)\times\R^d\ni t,\ba$ with $T>0$. In particular, one can choose smooth initial conditions for \eqref{variable}, or piecewise smooth with bounded variations and translation invariance so $d=1$ \cite{DafermosBook4,wagner-1987}. Then, the latter is equivalent to a well-defined Eulerian description 
for
\begin{equation}\label{variableshort}
(J^{-1},J^{-1}\bm,J^{-1}\bF,J^{-1}\rho_r^{r_\rho},J^{-1}K^{r_K},J^{-1}\bA^{-r_A}) 
\end{equation}
where we recall $\lambda\to\infty$, $G\equiv\frac{J \mu}{\lambda}\to0$ (no relaxation, no viscosity).
\end{proposition} 
To our knowledge, Prop.~\ref{prop:perfect} is the first well--posedness result for a genuine front-capturing modelling approach (based on a single PDE) to multi-material flows with possibly infinitely-many phases and sharp interfaces. Note that at this stage that it is natural to choose $r_A$ as a negative integer so that
\eqref{eq:Z_material_tensor_rA} simplifies to
\begin{equation}
\label{eq:A_material_tensor} \partial_t \bA^{r_A} = \left(\bA^{-r_A}(G\bF^{-1}\bF^{-T}\bA^{-r_A-1} + \dots
+ \bA^{-r_A-1}G\bA\bF^{-1}\bF^{-T})\bA^{-r_A}-\bA^{r_A}\right) \lambda^{-1}\,.
\end{equation}

Including viscosity in the multi-material flows above through a relaxation approach is possible as above when $\lambda>0$ is finite and $\mu\neq0\neq G$, but it has some limitations. Indeed, consider $(x,y)\to -x^\alpha\log y$ on $(\R_{>0})^2$: 
it is convex (jointly)
\begin{itemize}
 \item on $x>0$, $y\ge e^{\frac\alpha{1-\alpha}}>1$ with $0<\alpha<1$,
 \item on $x>0$, $e^{-1}<y\le e^{\frac\alpha{1-\alpha}}<1$ when $\alpha<0$, and
 \item on $x>0$, $0<y<e^{-1}$ when $\alpha>1$.
\end{itemize}
So $- G \log J$ 
cannot be jointly convex in $G^{r_G},J>0$ on one fixed domain with fixed $r_G$.
\begin{proposition}\label{prop:viscous}
Assume $\gamma>1$ and $r_K\in(0,\tfrac1{\gamma})$, $r_\rho>0$, $-r_A\in\NN_{>0}$. \\
The solutions to the system of quasilinear balance laws (\ref{eq:momentum_material_tensor}--\ref{eq:dilatation_material_tensor}--\ref{eq:deformation_material_tensor}--\ref{eq:reference_material_tensor}--\ref{eq:K_material_tensor}--\ref{eq:G_material_tensor}--\ref{eq:lambda_material_tensor}--\ref{eq:A_material_tensor}) that satisfy \eqref{piola_material} also satisfy one strictly convex extension on
\begin{equation}\label{variable2long2}
(e^{-\tfrac1{r_G-1}},\infty)\times\R^d\times\R^{d\times d}\times\R_{>0},\R_{>0},\R_{>0}^{d\times d},\R_{>0},\R_{>0} \ni
(J,\bm,\bF,\rho_r^{r_\rho},K^{r_K},\bA^{r_A},G^{r_G},\lambda^{r_\lambda}) 
\end{equation}
if 
$r_G>1$, another strictly convex extension on
\begin{equation}\label{variable2long1}
(e^{-1},e^{-\tfrac1{r_G-1}})\times\R^d\times\R^{d\times d}\times\R_{>0},\R_{>0},\R_{>0}^{d\times d},\R_{>0},\R_{>0} \ni
(J,\bm,\bF,\rho_r^{r_\rho},K^{r_K},\bA^{r_A},G^{r_G},\lambda^{r_\lambda}) 
\end{equation}
if 
$r_G<0$, and yet another convex extension on
\begin{equation}\label{variable2long3}
(0,e^{-1})\times\R^d\times\R^{d\times d}\times\R_{>0},\R_{>0},\R_{>0}^{d\times d},\R_{>0},\R_{>0} \ni
(J,\bm,\bF,\rho_r^{r_\rho},K^{r_K},\bA^{r_A},G^{r_G},\lambda^{r_\lambda}) 
\end{equation}
if 
$0<r_G<1$. \\
Therefore, they define an unequivocal Lagrangian description for \emph{viscoelastic} flows of immiscible infinitely-many materials insofar as they are the solutions to well-defined Cauchy problems on $(0,T)\times\R^d\ni t,\ba$ with $T>0$. \\
In particular, for one fixed $r_G$, one can choose smooth initial conditions for \eqref{variable} or piecewise smooth with bounded variations and translation invariance so $d=1$ \cite{DafermosBook4,wagner-1987},
\begin{itemize} 
 \item either with a uniform lower-bound $J_{1_+}>1$ on $J$ depending on $r_G>1$, 
 \item or with uniform lower-bound $J_{e^{-1}_+}>e^{-1}$ and upper-bound $J_{1_-}<1$ on $J$ depending on $r_G<0$,
 \item or with uniform upper-bound $J_{e^{-1}_-}<e^{-1}$ on  $J$ depending on $r_G\in(0,1)$.
\end{itemize}
The latter is still equivalent to a well-defined Eulerian description 
for
\begin{equation}\label{variableshort2}
(J^{-1},J^{-1}\bm,J^{-1}\bF,J^{-1}\rho_r^{r_\rho},J^{-1}K^{r_K},J^{-1}\bA^{r_A},J^{-1}G^{r_G},J^{-1}\lambda^{r_\lambda}) 
\end{equation}
however, the maximal time $T$ will now depend on $r_G$ (insofar as each symmetric-hyperbolic formulation breaks down 
when $J$ reaches one boundary value $e^{-1}$ or $1$ of its definition domain). Moreover, a formulation of the second principle of thermodynamics (i.e. a signed source term $H$ in the additional balance law for the ``energy'' extension) is furthermore guaranteed when the \emph{full} energy $E$ in \eqref{energy} is convex in \eqref{variableshort2} supplemented by the independent variable $\left(\log\frac{|A|}{|A|_-}\right)^\frac1\beta>0$, $|A|_-<|A|$ being a lower bound, 
i.e.
\begin{itemize}
 \item when $r_G<0$, using $\beta<0$, or
 \item when $0<r_G<1$, using $\beta\in\left(1-\tfrac1{r_G},0\right)$.
\end{itemize}
\end{proposition}
Note that in Prop.~\ref{prop:viscous}, we are not able to guarantee (an unequivocal solution with) a formulation of the second principle of thermodynamics in the case $J>1$ when one needs to consider $r_G>1$ for the joint convexity in $(G^{r_G},J)$ of the viscous potential term $-G\log J$, because the ``standard'' polyconvex approach to guarantee convexity of the other viscous potential term $- \frac{G}2 \log|A|$ in \eqref{energy} does not apply (constraints on $r_G$ contradict one another in the domain $J>1$).
One idea to overcome the previous limit is to double the variable $G$ into two different variables, as in the polyconvexity approach. We let this for future studies.
In any case, we have not (yet) considered adding a variable $\gamma-1>0$ (and a corresponding advection equation) to the system of balance laws~! Indeed, although natural e.g. in the multi-fluid context, it seems to us that requiring the joint convexity of $K\frac{J^{-(\gamma-1)}}{\gamma-1}$ in $K^{r_K},J,(\gamma-1)^{r_\gamma}$ would additionally impose excessive constraints on the model for the choice of $r_K,r_\gamma$.
That is, for the simulation of any multi-material flow by a general and genuine interface-capturing approach, difficulties remain in identifying a single system of PDEs 
endowed with a stable structure like symmetric-hyperbolic quasilinera PDEs typically issued from a hyperelastic approach, although a few cases can be mathematically grounded as we have show here.

\section{Conclusion}

A new interface-capturing approach for multi-material flows with sharp interfaces has been introduced.
It is based on a new system of PDEs to govern a viscoelastic 
mixture of possibly infinitely-many immiscible fluids and solids (phases) with Cauchy stress of neo-Hookean/Mooney-Rivlin (compressible, hyperelastic) type \eqref{formula}.
Viscous effects can be introduced through the relaxation of ``internal'' variables $\bZ$ as already introduced in \cite{Boyaval2021} for single-material viscoelastic flow (of Maxwell type).
Material heterogeneities are advected by the flow, as 
internal variables.

By contrast with other modelling approaches for immiscible fluid-solid mixtures, e.g. \cite{MR3481008,MR4038165}, our PDEs should deliver the promise of interface-capturing approaches and  genuinely define the multi-material flow without the help of additionnal hypotheses/tricks. 
Relying on the existing theory for quasilinear balance laws, we could prove that 
Cauchy problems (on the whole space $\R^3$) have well-defined solutions (smooth, or 1D)
\begin{itemize}
 \item on neglecting viscosity (like e.g. in \cite{MR3481008,MR4038165}),
 \item or for small variations of the dilatation/compression rate $J$ within a subdomain of $J>0$ where a symmetric-hyperbolic reformulation can be ensured.
\end{itemize}

So, to our knowledge, we have obtained here the first rigorous formulation of a genuine interface-capturing approach for immiscible multi-material flows. 
However, it is not obvious how to ensure the joint convexity of usual pressure (constitutive) laws with respect to both flow and material variables in general, even for the standard $\gamma$-law.
In particular, the present formulation of viscous effects cannot be guaranteed stable on the whole natural domain $J>0$ 
(the dilatation/compression rate should not vary too much within subdomains) and (a formulation of) the second principle of thermodynamics cannot be guaranteed in the dilatation case $J>1$ (with respect to the reference state).
Numerical tests could help confirm the stability limits in multi-material flows with viscous effects and then encourage the search for new formulations of the dissipative forces, or remain stable and then encourage for further mathematical tools to establish stability.


\end{document}